\documentclass[a4paper,12pt]{article}
\usepackage{hyperref} 
\usepackage[numbers]{natbib}

\usepackage{amsfonts,latexsym,amstext}
\usepackage{amsmath}
\usepackage{amssymb}
\usepackage[english]{babel}
\usepackage[latin1]{inputenc}
\usepackage{amsthm}
\usepackage{latexsym}
\usepackage{array}
\usepackage{mathrsfs}
\usepackage{amsxtra}
\usepackage{amscd}
\usepackage{mathtools}
\usepackage{verbatim}
\usepackage{enumerate}
\usepackage[a4paper,top=5cm,bottom=6cm,left=4cm,right=5cm,bindingoffset=7mm]{geometry}

\usepackage[usenames]{color}
\definecolor{beige}{rgb}{0.96, 0.96, 0.86}
\definecolor{airforceblue}{rgb}{0.36, 0.54, 0.66}
\definecolor{antiquefuchsia}{rgb}{0.57, 0.36, 0.51}
\definecolor{awesome}{rgb}{1.0, 0.13, 0.32}

\usepackage{xcolor}

\usepackage{fancyhdr}

\renewcommand{\tilde}{\widetilde}

\definecolor{red}{rgb}{1.0,0.0,0.0}

\definecolor{blu}{rgb}{0.0,0.0,1.0}

\definecolor{gre}{rgb}{0.03,0.50,0.03}

\definecolor{darkviolet}{rgb}{0.58, 0.0, 0.83}

\usepackage{enumerate}

\newtheorem{theorem}{Theorem}[section]

\newtheorem{proposition}[theorem]{Proposition}

\newtheorem{hypothesis}[theorem]{Hypothesis}

\numberwithin{equation}{section}

\def\qed{{\hfill\hbox{\enspace${ \square}$}} \smallskip}
\def\sqr#1#2{{\vcenter{\vbox{\hrule height .#2pt \hbox{\vrule
 width .#2pt height#1pt \kern#1pt \vrule
width .#2pt} \hrule height .#2pt}}}}
\def\square{\mathchoice\sqr54\sqr54\sqr{4.1}3\sqr{3.5}3}

\def\e{\varepsilon}

\def\ds{\begin{displaystyle}}
\def\eds{\end{displaystyle}}
\def\dis{\displaystyle }
\def\<{\left\langle }
\def\>{\right\rangle }

\def\dim{\noindent \hbox{{\bf Proof.} }}

\def\R{\mathbb R}
\def\N{\mathbb N}

\def\E{\mathbb E}
\def\P{\mathbb P}

\def\F{\mathcal F}

\def\calf{{\cal F}}

\def\calu{{\cal U}}

\def\to{\rightarrow}

\def\d{\mathtt d}

\title
{ Stochastic Maximum Principle for optimal advertising models  with delay and non-convex control space.}

\date{}

\makeatletter
\newcases{mycases}{\quad}{\hfil$\m@th\displaystyle{##}$}{$\m@th\displaystyle{##}$\hfil}{\lbrace}{.}
\makeatother

\author{{Giuseppina Guatteri}\\
Dipartimento di Matematica, Politecnico di Milano\\ via Bonardi 9, 20133\\ Milano, Italia\\
{giuseppina.guatteri@polimi.it}\\
\\
{Federica Masiero}
\\Dipartimento di Matematica e Applicazioni\\ Universit\`a di Milano-Bicocca\\ via Cozzi 55, 20125 Milano, Italia
\\{federica.masiero@unimib.it}}

\usepackage{cancel}

\begin{document}
\maketitle

\begin{abstract}
In this paper we study optimal advertising problems that models the introduction of a new product into the market in the presence of carryover effects of the advertisement and with memory effects in the level of goodwill. In particular, we let the dynamics of the product
  goodwill to depend on the past, and also on past 
  advertising efforts. We treat the problem by means of the stochastic Pontryagin maximum principle, that here is considered for a 
 class of problems where in the state equation either the state or the control depend on the past. Moreover the control acts on the martingale term, as in \cite{GuaMas}, but now the space of controls $ U $ can be chosen to be non-convex. The maximum principle is thus formulated using a first-order adjoint Backward Stochastic Differential Equations (BSDEs), which can be explicitly computed due to the specific characteristics of the model, and a second-order adjoint relation.

\end{abstract}



\section{Introduction}
We consider a stochastic model for problems of optimal advertising under uncertainty, so we have to study a stochastic version of advertising models. We start from the stochastic model introduced in \cite{GrossVisc} (see also \cite{GozMar} and \cite{GozMarSav}): we consider carryover effects of the advertising, which in the model reads as delay in the control, and with memory effects of the level of goodwill, which in the model reads as delay  in the state. We refer also to \cite{Hartl} for optimal advertising problems with memory effects in the level of goodwill. 

\noindent In our model the delay in the effect of advertising affects the martingale term of the state equation, 
\newline Namely we consider, following \cite{GrossVisc}, the controlled stochastic differential equation in $\R$ with  pointwise delay in the state and in the control
\begin{equation}
\left\{
\begin{array}
[c]{l}%
dx(t)  =\left[b_0u(t) - a_0x(t) -   {a_d} x(t-\d)\right]\,dt+\sigma_1 x(t)dW^1_t+\sigma_2 u(t-\mathtt{d})dW^2_t
, \\ \\
 x(\tau)  = x_0(\tau), \quad \tau \in [-\d,0]
.\\
\end{array}
\right.  
\end{equation}
Following the usual notations, see e.g. \cite{GrossVisc}, $x$ is the level of  goodwill, $a_0$ and $a_d$ are factors related to the  image de\-te\-rio\-ra\-tion in case of no advertising, $b_0$ is a constant representing an advertising effectiveness factor. The  diffusion term $\sigma_1  x(t)dW^1_t$ represents for the word of mouth communication,  with $\sigma_1\geq 0$ the so-called "advertising volatility"; while the second diffusion term keeps track of the delayed effect of the advertising effort $u$, and the constant $\sigma_2\geq 0$ is the effectiveness volatility of the communication. Besides, $x_0(0)$ is the initial goodwill level, while $x_0(\tau),\,\tau \in [-\d,0)$ is the history  of the goodwill level before the advertising campaign is started.

The functional to maximize, over all controls in $\mathcal {U}$, 
 is the following  profit, defined on finite horizon:
 \begin{equation}
\bar J(t,x,u)=\E \int_t^T \left(-c\left(u(s)\right) {+}l\left(x(s))\right)\right)\;ds +\E  \,  r(x(T)),
\end{equation}
where $c$ is the cost of advertising, $l$ is the current reward, and $r$ is the reward from the final goodwill.
Our purpose is to derive a maximum principle for such a problem  with non convex control space $U$, extending the results already  present in the literature for the convex case, see e.g. \cite{ChenWu}, where the stochastic maximum principle for control problems with pointwise delay in the state and in the control is studied,  see also \cite{GuaMas} and \cite{HuPe} where more general models are treated but the convexity of the control space $U$ is still required. We underline the fact  that stochastic control problems with diffusion depending on the control are difficult to treat; concerning problems related to  advertising we mention the recent paper \cite{deFeo}, where the author solves the problem by means of the dynamic programming approach, differently from here the case of pointwise delay cannot be treated, moreover it is taken into account a diffusion depending linearly on the control, but not on the delayed control. 
We stress that one novelty of this paper is to handle the new difficulty coming form the non convexity of $U$: also in the non delayed case the non convexity of $U$ makes the approach based on stochastic maximum principle more complicated. First we need to utilize the spike variation of the control and introduce the second variation to handle the control acting on the martingale term. We will therefore follow \cite{GuaMas2023} to introduce this additional step. Thanks to the specific characteristics of the state equation, we can derive a quasi-explicit form for both the first adjoint,
which results in an anticipating backward equation, see \cite{PengYang}, and the second adjoint, which is written using  the optimal cost and a simple auxiliary process.
Note that the specific case we address {is not considered in \cite{MeShiWaZha}} and does not completely fall within the hypotheses of \cite{GuaMas2023}, as we also consider the delay the control appearing in the martingale term. 

The structure of the paper is the following: in section  \ref{sec:opt-adv} we present the notation, the assumptions and the control problem.
Section \ref{sec:var} is devoted to collect the results on the first variation of the optimal state and also on the second variation, which we have to study since the space of controls is not convex. Section \ref{sec-duale}  concerns the first and second adjoint relations necessary to formulate the stochastic  maximum principle, which is stated and proved section \ref{sec max prin}.

\section{Stochastic control problem for delay equations appearing in advertising models}
\label{sec:opt-adv}
\subsection{Assumptions and preliminaries}\label{sec:AssPrel}

Let  $W(t)$ be a  $2$ dimensional brownian motion defined on a complete probability space $(\Omega, \calf, \P)$ and $(\mathcal{F}_t)_{t\geq 0}$ the natural filtration  generated by $W$, 
augmented in the usual way. We  fix a finite time horizon $T>0$  and  we set, for every $q \geq 1$ and for every $ 0\leq u\leq v\leq T$, the following spaces:
\begin{multline} L^0_{\mathcal{F}}(\Omega\times [u,v];\R^k)=: \{X:\Omega \times [u,v] \to  \R^k,  \text{  $(\mathcal{F}_t)$-progressive measurable\} } \\ \end{multline}
\begin{multline}   
L^q_{\mathcal{F}}(\Omega\times [u,v];\R^k) =: \\ \{X:\Omega \times [u,v] \to  \R^k,  \text{  $(\mathcal{F}_t)$-progr.meas., } \left( \E \int_u^v |X(t)|^q dt\right)^{1/q} < \infty \}\end{multline}
\begin{multline}   
L^q_{\mathcal{F}}(\Omega; C([u,v];\R^k)) =: \\ \{X:\Omega \times [u,v] \to  \R^k,  \text{  $(\mathcal{F}_t)$-progr.meas., } \Big( \E \sup_{t \in [u,v]} |X(t)|^q \Big)^{1/q} < \infty \}\end{multline}

 \subsection{Formulation of the problem}\label{sec-density}

We recall the  state equation  we are interested in:
\begin{equation}
\left\{
\begin{array}
[c]{l}%
dx(t)  =\left[b_0u(t) - a_0x(t) -   {a_d} x(t-\d)\right]\,dt+\sigma_1 x(t)dW^1_t+\sigma_2 u(t-\mathtt{d})dW^2_t
, \\ \\
x(\tau)  = x_0(\tau), \quad \tau \in [-\d,0),\;x(0)=x_0(0).\\
\end{array}
\right.  \label{eq-contr-rit}
\end{equation}

 We assume the following on the coefficients:
\begin{hypothesis}\label{ipotesibasic}
\begin{itemize}
\item[\textcolor{red}{(i)}] $a_0, b_0, \sigma_1,\, \sigma_2\in \R$;
\item[\textcolor{red}{(ii)}] the control strategy $u$, the advertisement in this case, belongs to the space{ $$\mathcal{U}:=\left\lbrace z\in  L^0_{\mathcal{F}}(\Omega\times [0,T], \R):z(t)\in U , \;\P-a.s.\right\rbrace $$}where $U$ is a bounded subset of $\R$ possibly non convex, in particular $U$ can be a bounded subset of $\N$ and this represents the realistic situation in which the advertisement is multiple of a given quantity;
\item[\textcolor{red}{(iii)}] $\d>0$ is the  delay with which the past of the state affects  the system;
 
\end{itemize}
\end{hypothesis}
We recall that the purpose is to maximize, over all controls in $\mathcal {U}$, the following profit, on finite horizon:
 \begin{equation}\label{costo-advertisement-maximize}
\bar J(t,x,u)=\E \int_t^T \left(-c\left(u(s)\right) {+}l\left(x(s))\right)\right)\;ds +\E  \,  r(x(T)),
\end{equation}

\begin{hypothesis}\label{ipotesibasic-costo}
 We assume the following:
\begin{itemize}
\item[(i)] according to the literature, see e.g. \cite{Tap}, $c:U\rightarrow \R$ is non linear, convex and locally lipschitz;
\item[(ii)] $l:[0,T]\times\R\to\R$ represents the current reward,  it is twice differentiable with at most linear growth;
\item[(iii)] $r:\R\to\R$ represents the foreseen reward from the final goodwill and it is twice differentiable and strictly increasing, with bounded second derivatives.
\end{itemize}
\end{hypothesis}
From now on, since we recall and apply the results in \cite{GuaMas2023} where a minimization problem is considered, we focus on the problem (equivalent to the original one of maximizing the profit $\bar J$ given in \eqref{costo-advertisement-maximize} since $J=-\bar J$) of minimizing, over all admissible control in $\calu$, the cost functional
\begin{equation}\label{costo-advertisement}
 J(t,x,u)=\E \int_t^T \left(c(u(s)) -l(x(s))\right)ds -\E  \,r(x(T)),
\end{equation}
We are going to formulate necessary conditions for optimality.
\section{First and second order variations of the optimal state equation} \label{sec:var}
Since $U$ is not necessarily convex we  use  the spike variation method, see \cite{Peng}, see also \cite{GuaMas2023} for the case of delay equation.
\newline Let  $(u,x)$ be an optimal couple: that is  ${u}$ is an optimal control and ${x}$  its corresponding optimal trajectory, that is  solution to equation \eqref{eq-contr-rit}.  The {\em spike} variation works as follows:
for $\e >0$  we consider an interval $V_\e\subset [0,T]$, with $m(V_\e)= \e$, where $m$ is the Lebesgue measure.
Let $v \in U$ 
we set
\begin{equation}\label{spike_def}
u^\e(t) = \left\{ \begin{array}{ll}
{u} (t)  & t \in [0,T] \setminus V_\e\\
v & t \in V_\e
\end{array}\right.
\end{equation}
We are going to derive a maximum principle for the control problem with state equation \eqref{eq-contr-rit} and cost functional \eqref{costo-advertisement} in the case of non convex space of controls.
We will denote by $x^\e$ the solution of \eqref{eq-contr-rit} corresponding to $u^\e$ and $\delta u$ the spike variation of the control, i.e. $\delta u(t):= (u(t) -v))I_{V_\e}(t)$.
\newline We can write the equation for the first order variation of the state:
\begin{equation}\label{var.prima.fin.dim}
\left\{ \begin{array}{ll}
dy^\e(t) =  \left[-a_0 y^\e(t) - a_dy^\e(t-\d)\right]\,dt+\sigma_1 y^\e(t)dW^1_t+\sigma_2 \delta u(t-\d)dW^2_t
\quad    t \in [0,T], \\ \\
  y^\e(0)=0, \quad y^\e(\tau)=0, \quad-\mathtt{d} \leq \tau <0.
\end{array} \right.
\end{equation}
where $\delta u^\e(t-\d)= (u(t-\d) -v)I_{V_\e} (t-\d) $
and the second variation is, for $ t \in [0,T]$,
\begin{equation}\label{var.seconda.fin.dim}
\left\{ \begin{array}{ll}
dz^\e(t) =    \left[-a_0 z^\e(t) - a_dz^\e(t-\d)\right]\,dt+b_0\delta u(t)dt+\sigma_1 z^\e(t)dW^1_t\\  \\
  z^\e(0)=0, \quad z^\e(\tau)=0, \quad-\mathtt{d} \leq \tau <0.
\end{array} \right.
\end{equation}
It is well known  that such equations are well posed in  $L^2_{\mathcal{F}}(\Omega; C([0,T];\R))$ (see e.g. \cite{Moha} for a general theory on stochastic delay equations).

The following asymptotic behaviors hold, see \cite[Theorem 3.3]{GuaMas2023}  for the delayed system and the classic result in \cite[Theorem 4.4]{YongZhou},
\begin{equation}\label{eq:stima_varI-concreto}
\E \sup_{t \in [0,T]} |y^\e (t)|^{2} = O(\e), 
\qquad
\E \sup_{t \in [0,T]} |z^\e (t)|^{2} = O(\e), 
\end{equation}
\begin{equation}\label{eq:approx-stateI-concreto}
\E \sup_{t \in [0,T]} |x^\e(t)-x(t)-y^\e (t)|^{2} = O(\e^{2}), 
\end{equation}
\begin{equation}\label{eq:approx-stateII-concreto}
\E \sup_{t \in [0,T]} |x^\e(t)-x(t)-y^\e (t)-z^\e(t)|^{2} = o(\e^{2}), 
\end{equation} 
Moreover,  $\forall\,p\geq 1$, 
\begin{align}\label{stima sup y}
        & \E \sup_{t\in [-\mathtt{d},T]} |y^{\e}(t)|^p < +\infty,\qquad\E \sup_{t\in [-\mathtt{d},T]} |z^{\e}(t)| ^p < +\infty, 
      \end{align} 
\section{ Maximum Principle for the stochastic delay equation}\label{sec-duale}

In this Section, first we formulate the first order adjoint equation in \ref{subsec-dualeI-anticp}, then we pass to the second order adjoint in \ref{sec-dualeII}, and finally we formulate the stochastic maximum principle in \ref{sec max prin}.
\subsection{First order duality relation}\label{subsec-dualeI-anticp}
 Following \cite{GuaMas2023} we define the first order adjoint equation of \eqref{var.prima.fin.dim}, that is the anticipated backward equation 
 \begin{equation}\label{aggiunta-advertisement}
  \left\lbrace\begin{array}{l}
p(t)= r_x(x(T)) + \dis\int_t^T l_x\left( x(s)\right)\, ds
- a_0 \dis\int_t^Tp(s)\,ds- a_d\dis\int_t^T\E^{\calf_s}p(s+\d)\, ds\\ \\ \qquad+ \sigma_1\dis\int_t^T \E^{\calf_s}q^{1}(s+\d)\,ds+\dis\int_t^Tq^{1}(s)dW^1_s+\dis\int_t^Tq^{2}(s)dW^2_s\\ \\
p(T-\tau)=0, \; \text{for all }\tau  \in [-\mathtt{d},0 ), \\ \\  q^1(T-\tau)= q^2(T-\tau)= 0 \quad \text{a.e. }\tau  \in [-\d,0).
 \end{array}
 \right.
 \end{equation}
 that admits a unique solution $(p,q)\in
L^2_{\mathcal{F}}(\Omega; C([0,T];\R)) \times
 L^2_{\mathcal{F}}(\Omega\times [0,T];\R^2)$, see e.g. \cite{PengYang} where we write $q=\left(q^{1},q^{2} \right)$.

Since  equation \eqref{aggiunta-advertisement} is linear, in view of the future times conditions we have an explicit, recursive, formulation for the solution $(p,q)$:
\begin{proposition}
    Assume  Hypotheses \ref{ipotesibasic} and \ref{ipotesibasic-costo}. Define for every  $\ k= 0, \dots, \Big[\frac{T}{\d}\Big]$,
    \[p^k(t):= p_{|_{[(T-(k+1)\d) \vee 0, (T-k\d)  ]}(t)} \]\text{ and } \[q^{k,1}(t):= q^1_{|_{[(T-(k+1)\d) \vee 0,(T-k\d) ]}}(t), q^{k,2}(t):= q^2_{|_{[(T-(k+1)\d) \vee 0,(T-k\d)  ]}}(t).\]
    Then $(p^k,q^k)$ solve:
 \begin{equation}\label{aggiunta-advertisement- ricur}
  \left\lbrace\begin{array}{ll}
p^0(t)=   e^{- a_0 (T -t)} \E^{\F_t}(r_x(x(T))) +   \E^{\F_t} \dis\int_t^T  e^{- a_0 (s -t)} l_x\left( x(s)\right)\, ds,  \quad t \in [T-\d,T], \\  \\ q^{0,1}(t)=   e^{- a_0 (T -t)}L_1(t)  - \dis\int_t^T   e^{- a_0 (s -t)} K_1(s,t)\, ds, \qquad \quad\text{      for a.e. } t \in [T-\d,T],  \\
 q^{0,2}(t)=   e^{- a_0 (T -t)}L_2(t)   - \dis\int_t^T   e^{- a_0 (s -t)} K_2(s,t)\, ds,\qquad \quad\text{      for a.e. } t \in [T-\d,T]. \\
p^k(t)=  e^{- a_0 (T-k \d -t)}
p^{k-1}(T-k \d)- a_d \dis\int_t^{T-k \d} 
e^{-a_0 (s-t)} \E^{\calf_s} p^{k-1}(s+\d)\, 
ds 
\\ \\   + \sigma_1\dis\int_t^T  e^{-a_0 (s-t)}  \E^{\calf_s}q^{k-1,2}(s+\d)\,ds+\dis\int_t^T   e^{-a_0 (s-t)} q^{k,1}(s)dW^1_s \\  \\ +   \dis\int_t^T    e^{-a_0 (s-t)} q^{k,2}(s)dW^2_s, \quad t \in [(T-(k+1)\d) \vee 0,(T-k\d) ],  \, k:1, \dots , \Big[\frac{T}{\d}\Big] \\ 
 \end{array}
 \right.
 \end{equation}
 Where $L_1$ and $L_2$ are given by:
 \begin{equation}\label{rappresentazione}
    \E^{\calf_t} ( r_x(X(T)))=  \E ( r_x(X(T))) + \int_t^T L_1(s) \, d W^1(s)  + \int_t^T L_2(s) \, d W^1(s). \end{equation} 
   and $K_1$ and $K_2$ are given by:
    \begin{equation}\label{rappresentazione costo corrente}
    \E^{\calf_t} ( l_x(s))=  \E ( l_x(s)) + \int_{T-\d}^t K_1(s,\tau) \, d W^1(\tau)  + \int_{T-\d}^t K_2(s,\tau) \, d W^2(\tau) . \end{equation} 
\end{proposition}
\dim Notice that, to avoid trivialities, we are taking  $\d < T$.
\newline Let us consider the first order adjoint equation \eqref{aggiunta-advertisement}, (, , for $t\in (T-\d,T]$: in view of the fact that $p(t)=q^i(t)=0,\, i=1,2$ for a.e. $t\in [T,T+\d]$, for $t\in [T-\d,T]$ it can be rewritten as a linear BSDE 
\begin{equation}\label{aggiunta-advertisement-1}
  \left\lbrace\begin{array}{l}
p(t)= r_x(x(T)) + \dis\int_t^T l_x\left( x(s)\right)\, ds
- a_0 \dis\int_t^Tp(s)\,ds+\dis\int_t^Tq^{1}(s)dW^1_s\\ \\  \qquad+\dis\int_t^Tq^{2}(s)dW^2_s\\ \\
p(T-\tau)=0, \; \text{for all }\tau  \in [-\mathtt{d},0 ), \\ \\  q^1(T-\tau)= q^2(T-\tau)= 0 \quad \text{a.e. }\tau  \in [-\d,0),
 \end{array}
 \right.
 \end{equation}
 and its solution for $t \in [T-\d,T],$  is given by
 \begin{align}
     \label{aggiunta-advertisement- ricur-1}
p^0(t) &=  e^{- a_0 (T -t)} r_x(x(T)) +\dis\int_t^T   e^{-a_0 (s-t)} l_x(x(s))\, ds 
+\dis\int_t^T   e^{-a_0 (s-t)} q^{0,1}(s)dW^1_s  \nonumber \\ 
&+   \dis\int_t^T    e^{-a_0 (s-t)} {q^{0,2}}(s)dW^2_s
 \end{align}
 Formula \eqref{aggiunta-advertisement- ricur} is just the rearrangement of the variation of constant formula, while  \eqref{rappresentazione}  and \eqref{rappresentazione costo corrente} is an application of the Martingale Representation Theorem, see also \cite{HuPeng}.
\qed

Using \eqref{var.prima.fin.dim}, \eqref{var.prima.fin.dim} and \eqref{aggiunta-advertisement}, we deduce the first duality relation:
\begin{proposition}\label{prop:Ito}
Assume \ref{ipotesibasic} and \ref{ipotesibasic-costo} then: 
       \begin{align}\label{duality relation y}
\E p(T)y^{\e}(T))=\E r_x(x(T))y^{\e}(T) = &-\E\dis\int_0^T[l_x(x(s))y^\e(s)+q^{2}(s)\delta u(s-\d)]\,ds
\end{align}\
     \begin{align}\label{duality relation z}
&\E p(T)z^{\e}(T))=\E r_x(x(T)) z^\e(T)=-\E\dis\int_0^Tl_x(x(s))z^\e(s)\,ds + \E\dis\int_0^Tp(s)\delta u(s)\,ds  
\end{align} 
Moreover the cost can be written as follows:
 \begin{align}\label{costo-espansione-I-gen}
& J(u)-J(u^\e) = 
\E\int_0^T \left(-
\delta l (t)+  b_0 p(t) \delta u(t) - \sigma_2 q^{2}(t) \delta u (t-\d)\right) \,dt \\ 
 &\nonumber +\frac{1}{2} \E \, r_{xx}(x(T))(y^\e(T))^2 
+\frac{1}{2}\E \int_0^T l_{xx}(t,x(t))(y^\e(t))^2\, dt 
\end{align}
\end{proposition}
\dim  See \cite[Propositions 4.1 and 4.2]{GuaMas2023}.
\qed
\subsection{Second order duality relation} \label{sec-dualeII}
Following  \cite[Theorem 4.9]{GuaMas2023} we introduce the process $P$,  called $P_{0,0}$ in \cite{GuaMas2023}, that in this simplified case can be explicitely defined by means of a auxiliary process $\tilde{y}^{s,1}$ (here $s$ indicates the initial time and $1$ is the initial condition) that solves the following equation
\begin{equation}\label{y.dato}
\left\{ \begin{array}{ll}
d{\tilde y}^{s,1}(t) = (\displaystyle -a_0{\tilde y}^{s,1}(t) - a_d\tilde y^{s,1} (t- \d)) \, dt    +\displaystyle\sigma_1\tilde{y}^{s,1}(t) d W(t), \quad    t \in [s,T], \\ \\
  {\tilde y}^{s,1}(s)=1, \quad y(\tau)=0, \quad -\d \leq \tau <s.
\end{array} \right.
\end{equation}
 By means of the solution process $\tilde y^{s,1} \in  L^2_{\mathcal{F}}(\Omega; C([0,T];\R)) $, the process  $P$ can be now defined by  
\begin{align}    
\label{teoP00formula}
     & P(s)=  \E ^{\F_s} r_{xx}(x(T))\tilde y^{s,1}(T)^2
 + \E ^{\F_s}\displaystyle \int_s^T  l_{xx}(t)\tilde y^{s,1}(t)^{2}\, dt.
\end{align}

Also in this case, the solution $\tilde y^{1,s}$  of equation \eqref{y.dato} can be explicitly defined by recursion:
\begin{proposition}
    Assume  Hypotheses \ref{ipotesibasic} and \ref{ipotesibasic-costo}. Define
    \[\tilde y^k(t):= \tilde y^{s,1}_{|_{[(s+ k\d) ,(s+ (k+1)\d)\wedge T]}}(t), \text { for every } k: 0, \dots, \Bigg[\frac{T}{\d}\Bigg]\]
    Then $\tilde y^k$ solves
\begin{equation}\label{stato- ricur}
  \left\lbrace\begin{array}{ll}
  \tilde y^0(t)= e^{(a_o - \frac{\sigma_1^2}{2}(t-s) +\sigma_1 (W^1(t)-W^1(s))} \quad t \in [s,(s+\d) \wedge T] \\ \\
\tilde y^k(t)=  e^{ -a_0 (t- (s+ k\d))}\tilde  y^{k-1}( s+ k\d) - a_d \dis\int^t_{ (s+ k\d)}  e^{-a_0 (t-r)}\tilde  y^{k-1}(r-\d)\, ds  \\ 
+ \sigma_1\dis\int^t_{ (s+ k\d)}   e^{a_0 (t-r)}\tilde  y^{k}(s)dW^1_s  \quad t \in [(s+k\d), (s+(k+1)\d)\wedge T], \,  k= 1, \dots, \Big[\frac{T}{\d}\Big]
 \end{array}
 \right.
 \end{equation}
 \end{proposition}
 \dim
Thanks to the initial condition in \eqref{y.dato}, we have that for every $t \in [s,s+\d \wedge T]$:
 \begin{equation}
  \left\lbrace\begin{array}{l}
d \tilde{y}^{s,1}(t)=   - a_0 \tilde{y}^{s,1} (t) \, dt + \sigma_1 \tilde{y}^{s,1}(t)dW^1_t \\
\tilde{y}^{s,1}(s)=1,
 \end{array}
 \right.
 \end{equation}
 Therefore $\tilde y^{0}(t)=  \tilde{y}^{s,1}_{|_{[s ,(s+\d) 
 \wedge T]}}(t)= 
 e^{(a_o - \frac{\sigma_1^2}{2}(t-s) +\sigma_1 (W^1(t)-W^1(s))}$.
 Then the general  case for $t> s+ \d$ follows by the application of the constant variation formula on any interval $ [(s+k\d), (s+(k+1)\d)\wedge T]$.
 \qed

The expansion of the cost then becomes:
\begin{proposition}\label{prop:Ito 2}
Assume \ref{ipotesibasic} and \ref{ipotesibasic-costo}. Then following expansion holds true: 
 \begin{align}\label{costo-espansione-II-gen}
& J(u)-J(u^\e) = 
\E\int_0^T \left(-
\delta c (s)+ b_0 p(s)\delta u(s) - \sigma_2 q^{2}(s) \delta u (s-\d)\right) (t)\,ds \\ 
 &\nonumber +\frac{1}{2}  \E\int_0^T  \sigma_2^2( \delta u(s-\d))^2  P(s) \,ds + o(\e) 
\end{align}
\end{proposition}

\dim
The proof is in \cite{GuaMas2023} [Proposition 4.9, Theorem 4.10, Theorem 4.11] and is a consequence of an approximation procedure based on the regularization of the dirac measure $\delta_{\d}$ and the well known reformulation \textcolor{red}{of} the delayed problem in infinite dimensions.
\qed

\subsection{The Stochastic Maximum principle}
\label{sec max prin}
We are now in the position to state our main result.
\begin{theorem}\label{teoMaxPrinc-gen}
Under assumptions \ref{ipotesibasic} and \ref{ipotesibasic-costo}, any optimal couple  $(x,u)$ satisfies 
the following variational inequality $ $:
 \begin{multline}\label{principio max formula}
 b_0(u(t)-v) p(t)
 - \sigma_2  (u(t)-v)\E^{\calf_t}q^{2}(t+\d) 
-(c(u(t))-c(v)) \\  \\+ \frac{1}{2} \sigma_2   (u(t)-v)^2 \E^{\calf_t}P(t+\d)  \leq 0, {\qquad\forall\,v\in U,\  a.e.\, \P\text{-a.s.}}
  \end{multline}
  \medskip
  
Where $ (p,q) \in L^2_{\mathcal{F}}(\Omega; C([0,T];\R)) \times
 L^2_{\mathcal{F}}(\Omega\times [0,T];\R^2)$, is the solution of the first order adjoint equation \eqref{aggiunta-advertisement},  and $P$ is the second adjoint process given by \eqref{teoP00formula}.
\end{theorem}
\dim
Let $t \in [0,T]$, and $V_\e=[t, t+\e]$, then
 \begin{align}
 &J(u)-J(u^\e)  = 
\E\int_0^T \left(-
\delta c (s)+ b_0 p(s)\delta u(s) - \sigma_2 q^{2}(s) \delta u (s-\d)\right) (t)\,ds \\ 
 &\nonumber +\frac{1}{2}  \E\int_0^T  \sigma_2^2( \delta u(s-\d))^2  P(s) \,ds= \\\nonumber
& \E\int_0^T \left(-[
  (c(u(s))-v) + b_0 p(s)(u(s) -v)] I_{V_\e}(s)-\sigma_2 q^{2}(s) (u (s-\d) -v) I_{V_\e}(s-\d)\right)\,ds \\& \nonumber+\frac{1}{2}  \E\int_0^T  \sigma_2^2( u(s-\d) -v)^2  P(s) I_{V_\e}(s-\d)\,ds= \\&\nonumber -
 \E\int_t^{t+\e} \left(
  (c(u(s))-v) I_{V_\e}(s)- b_0 p(s)(u(s) -v)\right)\,ds \\& \nonumber - \E\int_{t}^{t+\e}  \E^{\calf_{t}}\sigma_2 q^{2}(s'+\d) (u (s') -v) \, ds  + \frac{1}{2} \E\int_{t}^{t+\e} \E^{\calf_{t}}\sigma_2^2( u(s') -v)^2  P(s'+\d)]\,ds'.  
\end{align}

{Letting $\e$ tends to $0$, by standard arguments we deduce \eqref{principio max formula}}.

\bigskip

\textbf{Acknowledgements}: {The authors are grateful to 
the GNAMPA (Gruppo Nazionale per l'Analisi Matematica,
la Probabilit\`a e le loro Applicazioni) for the financial support.}


   \appendix

\end{document}